\documentclass[smallextended]{svjour3}
\usepackage{amsfonts}

\usepackage{latexsym}
\usepackage{amssymb}
\usepackage{amsmath}
\usepackage{graphicx}
\usepackage{hyperref}
\usepackage{caption}
\usepackage{subcaption}

\renewcommand{\qed}{\hfill{\ \ \rule{2mm}{2mm}} \vspace{0.2in}}

\newcommand{\ind}{1\hspace{-2.3mm}{1}}

\renewcommand{\thefigure}{\arabic{figure}}
\begin{document}

\title{Robust Domination in Random Graphs}
\author{ \textbf{Ghurumuruhan Ganesan}
\thanks{E-Mail: \texttt{gganesan82@gmail.com} } \\
\ \\
IISER Bhopal}
\date{}
\maketitle

\begin{abstract}
In this paper, we study ``robust" dominating sets of random graphs that retain the domination property even if a  small \emph{deterministic} set of edges are removed. We motivate our study by illustrating with examples from wireless networks in harsh environments. We then use the probabilistic method and martingale difference techniques to determine sufficient conditions for the asymptotic optimality of the robust domination number. We also discuss robust domination in sparse random graphs where the number of edges grows at most linearly in the number of vertices.

\vspace{0.1in} \noindent \textbf{Key words:} Random Graphs; Robust Domination; Sparse Regime.

\vspace{0.1in} \noindent \textbf{AMS 2000 Subject Classification:} Primary: 06C05; 05C62;
\end{abstract}

\bigskip

\setcounter{equation}{0}
\renewcommand\theequation{\thesection.\arabic{equation}}
\section{Introduction} \label{sec_intro}
Domination of graphs is an important topic from both theoretical and application perspectives and has been extensively studied in the random graph context as well. Throughout, random graphs refer to the Bernoulli or Erd\"os-R\'enyi random graph~\(G\) obtained by allowing each edge in the complete graph on~\(n\) vertices to be present with a certain probability~\(p,\) independent of the other edges (for formal definitions, please refer to Section~\ref{sec_dom}).  In~\cite{wieland}, two point concentration for the domination number of~\(G\) is obtained for the case when~\(p\) is essentially a constant and this concentration phenomenon was extended for a wide range of~\(p\) in~\cite{glebov}. Since then many other variants of domination have also been studied (see for e.g.~\cite{clark}\cite{wang}). 


Dominating sets also occur naturally in the design of wireless networks. In~\cite{wu}, an early application of dominating sets is explored for routing in ad hoc wireless networks devoid of any central control. The nodes belonging to the dominating sets are interpreted as ``gateways" through which any two nodes in the network can communicate with minimal delay. This was extended to higher dimensional wireless networks in~\cite{zou} and for a survey on the usage of domination in communications, we refer to~\cite{du}.

Ad hoc networks are especially fragile in terms of linkage in the sense that the sensors are continually moving around and so links may break or form randomly. Moreover, due to environmental constraints like shadowing and fading, it may happen that links between certain nodes are simply not feasible. In such a situation, it is of natural interest to know whether the domination property is still retained and this is the topic of study in this paper. We consider dominating sets in the Bernoulli random graph~\(G\) and are interested in obtaining ``robust" dominating sets that retain the domination property even if a small deterministic set of edges are removed.  We obtain sufficient conditions for asymptotic optimality of the robust domination number in terms of the maximum vertex degree and the number of edges of the graph that has been removed from~\(K_n.\)

The paper is organized as follows: In Section~\ref{sec_dom}, we state and prove our main result regarding robust domination number in random graphs for the dense regime. We use the probabilistic method to establish sufficient conditions for asymptotic optimality. Next, in Section~\ref{sec_sparse}, we also discuss our results for the robust domination in the sparse regime.


\setcounter{equation}{0}
\renewcommand\theequation{\thesection.\arabic{equation}}
\section{Robust Domination} \label{sec_dom}
Let~\(K_n\) be the complete graph on~\(n\) vertices and let~\(\{Z(f)\}_{f \in K_n}\) be independent random variables indexed by the edge set of~\(K_n\) and with distribution
\begin{equation}\label{x_dist}
\mathbb{P}(Z(f) = 1) = p = 1-\mathbb{P}(Z(f) = 0)
\end{equation}
where~\(0 < p < 1.\)  Let~\(G\) be the random graph formed by the union of all edges~\(f\) satisfying~\(Z(f) = 1\) and let~\(H\) be any \emph{deterministic} subgraph of~\(K_n\) with~\(m = m(n)\) edges and a maximum vertex degree of~\(\Delta = \Delta(n).\)

A set~\({\cal S} \subset V\) is said to be a \emph{dominating set} of~\(G \setminus H\) if each vertex in~\(V \setminus {\cal S}\) is adjacent to at least one vertex in~\({\cal S}\) in the graph~\(G \setminus H.\) We also say that~\({\cal S}\) is a~\(H-\)\emph{robust} dominating set or simply a robust dominating set. The~\(H-\)\emph{robust domination number} or simply the robust domination number is defined to be the minimum size of a dominating set in~\(G \setminus H\) and is denoted by~\(\gamma(G \setminus H).\)

We seek conditions on~\(H\) so that the robust domination number~\(\Gamma_n := \gamma(G\setminus H)\) and the actual domination number~\(\gamma(G) \leq \Gamma_n\) are of the same order. Intuitively, if~\(H\) is sparse, then we expect~\(\Gamma_n\) and~\(\gamma(G)\) to be close with high probability, i.e., with probability converging to one as~\(n \rightarrow \infty.\) This is illustrated in our first result below that obtains bounds for~\(\Gamma_n\) in terms of the maximum vertex degree~\(\Delta\) of the graph~\(H.\) For~\(0 < x,y < 1\) we define~\(u_n(x,y):= \frac{\log(nx)}{|\log(1-y)|}\) and set~\(u_n := u_n(p,p).\) Moreover, we use the notation~\(a_n = o(b_n)\) to denote that~\(\frac{a_n}{b_n} \longrightarrow 0,\) as~\(n \rightarrow \infty.\)
\begin{lemma}\label{thm_inter} The following properties hold:\\
\((a)\) Let~\(\lambda_a := np\) and~\(\lambda_b := n|\log(1-p)| > \lambda_a.\) For every~\(\theta > 2,\) there exists a~\(\lambda_0  = \lambda_0(\theta) > 0\) such that if~\(\frac{\lambda_0}{n} \leq p \leq 1-\frac{1}{n^3},\) then
\begin{equation}\label{dom_bound_low}
\mathbb{P}\left(\Gamma_n \geq u_n \left(1-\frac{\theta\log\log{\lambda_b}}{\log{\lambda_a}}\right)\right) \geq 1-\exp\left(-\frac{3n}{8} \cdot \frac{(\log{\lambda}_b)^{\theta}}{\lambda_a}\right).
\end{equation}
\((b)\) Let~\(\Delta\) be the maximum vertex degree of~\(H\) and suppose~\(np \longrightarrow \infty\) and\\\(p \leq p_0\) for some constant~\(0 < p_0 <1.\) For every~\(\epsilon > 0\) and all~\(n\) large,
\begin{equation}\label{dom_bound_up_ax}
\mathbb{P}\left(\Gamma_n \leq u_n(1+\epsilon) + \Delta\right) \geq 1- \frac{1}{\log(np)}.
\end{equation}
Consequently if~\(\Delta = o(u_n)\) and~\(p \leq p_0\) is such that~\(np \longrightarrow \infty,\) then~\(\frac{\Gamma_n}{u_n} \longrightarrow 1\) in probability as~\(n \rightarrow \infty.\)
\end{lemma}
The condition~\(np \longrightarrow \infty\) ensures that~\(G\) is reasonably dense in terms of its vertex degrees. In particular if the edge probability~\(p\) is a constant, then~\(u_n\) is of the order of~\(\log{n}\) and moreover both~\(\lambda_a\) and~\(\lambda_b\) are of the order of~\(n.\) Setting~\(\theta = 3\) in~(\ref{dom_bound_low}), we then get that~\(\Gamma_n \geq u_n \left(1- O\left(\frac{\log\log{n}}{\log{n}}\right)\right)\) with probability at leat~\(1-e^{-C(\log{n})^3}\) for some constant~\(C > 0.\) Similarly, the final statement~(\ref{dom_bound_up_ax}) implies that if~\(\Delta = o(\log{n}),\) then~\(\Gamma_n \leq u_n(1+2\epsilon)\) with high probability, for any arbitrary constant~\(\epsilon > 0.\)

Combining the observations of the previous paragraph, we get that~\(\Gamma_n \sim u_n\) with high probability, where we use the notation~\(a_n \sim b_n\) to denote that\\\(\frac{a_n}{b_n} \longrightarrow 1\) as~\(n \rightarrow \infty.\) Thus the robust domination number satisfies\\\(\Gamma_n \sim u_n \sim \gamma(G)\) and is therefore asymptotically equal to the ``ideal" domination number~\(\gamma(G),\) with high probability.




\emph{Proof of Lemma~\ref{thm_inter}~\((a)\)}: Since~\(\gamma(G \setminus H) \geq \gamma(G),\) it suffices to lower bound~\(\gamma(G)\) and for completeness, we give a small proof using a union bound argument covering all possibilities, as in~\cite{wieland}~\cite{glebov}. Specifically, letting~\(\lambda_a := np, \lambda_b = n|\log(1-p)|\)  and~\(t_n := \frac{\log{\lambda_a} - \theta\log\log{\lambda_b}}{|\log(1-p)|}\) vertices, we show that there exists a dominating set containing~\(t_n\) vertices, with high probability.

We begin with upper and lower bounds for~\(t_n.\) Using~\(|\log(1-p)| > p,\) we see that~\(t_n \leq \frac{\log{\lambda_a}}{p} = n\frac{\log{\lambda_a}}{\lambda_a} < \frac{n}{4}\) if~\(\lambda_a \geq \lambda_0,\) a sufficiently large absolute constant. Moreover, we have that~\(\lambda_a > (\log{\lambda_b})^{2\theta}\) for all~\(n \geq N_0 = N_0(\theta)\) large, provided~\(np \geq \lambda_0 = \lambda_0(\theta)\) is large. Indeed if~\(p \leq \frac{1}{2},\) then using~\(|\log(1-p)| < 2p,\) we get that~\(\lambda_a  - (\log{\lambda_b})^{2\theta} \geq np-(\log(2np))^{2\theta} >0\) if~\(np \geq \lambda_0 = \lambda_0(\theta)\) is sufficiently large. On the other hand if~\(\frac{1}{2} \leq p \leq 1-\frac{1}{n^3},\) then~\(\lambda_a = np > \frac{n}{2}\) and~\[ (\log{\lambda_b})^{2\theta} =\left(\log n + \log\log\left(\frac{1}{1-p}\right)\right)^{2\theta} \leq (\log{n} + 6\log\log{n})^{2\theta} < \frac{n}{2}\] for all~\(n \geq N_0 = N_0(\theta)\) large. Summarizing, we get that
\begin{equation}\label{t_bounds}
\frac{\log{\lambda_a}}{2|\log(1-p)|} < t_n < n\frac{\log{\lambda_a}}{\lambda_a} < \frac{n}{4}
\end{equation}
for all~\(n\) large.

Let~\({\cal S}\) be any set containing~\(t_n\) vertices. For a vertex~\(v \in  {\cal S}^c,\) the probability that~\(v\) is not adjacent to any vertex of~\({\cal S}\) is~\((1-p)^{t_n} = \frac{(\log{\lambda_b})^{\theta}}{\lambda_a}.\) Thus the vertex~\(v\) is adjacent to some vertex of~\({\cal S}\) with probability~\(1-\frac{(\log{\lambda_b})^{\theta}}{\lambda_a}\) and so if~\(E_{dom}({\cal S})\) is the event that~\({\cal S}\) is a dominating set, then using the fact that the complement set~\({\cal S}^c\) has~\(n-t_n \geq \frac{3n}{4}\) vertices (see~(\ref{t_bounds})), we get that
\[\mathbb{P}\left(E_{dom}\left({\cal S}\right)\right) \leq \left(1-\frac{(\log{\lambda_b})^{\theta}}{\lambda_a}\right)^{\frac{3n}{4}} \leq \exp\left(-\frac{3n}{4} \cdot \frac{(\log{\lambda_b})^{\theta}}{\lambda_a}\right).\]

Since there are~\({n \choose t_n} \leq \left(\frac{ne}{t_n}\right)^{t_n} = \exp\left(t_n  \log\left(\frac{ne}{t_n}\right)\right)\) sets of size~\(t_n,\) we use the union bound and the bounds in~(\ref{t_bounds}), to see that the probability that there exists a dominating set of size at most~\(t_n\) is bounded above by
\begin{eqnarray}
&&\exp\left(t_n  \log\left(\frac{ne}{t_n}\right)\right)\exp\left(-\frac{3n}{4} \cdot \frac{(\log{\lambda_b})^{\theta}}{\lambda_a}\right) \nonumber\\
&&\;\;\leq \;\;\exp\left(n\frac{\log{\lambda_a}}{\lambda_a} \log\left(\frac{2ne|\log(1-p)|}{\log{\lambda_a}}\right)\right)\exp\left(-\frac{3n}{4} \cdot \frac{(\log{\lambda_b})^{\theta}}{\lambda_a}\right) \nonumber\\
&&\;\;=\;\;\exp\left(n\frac{\log{\lambda_a}}{\lambda_a} \log\left(\frac{2e\lambda_b}{\log{\lambda_a}}\right)\right)\exp\left(-\frac{3n}{4} \cdot \frac{(\log{\lambda_b})^{\theta}}{\lambda_a}\right) \nonumber\\
&&\;\;\leq\;\;\exp\left(n\frac{\log{\lambda_a}}{\lambda_a} \log\left(2e\lambda_b\right)\right)\exp\left(-\frac{3n}{4} \cdot \frac{(\log{\lambda_b})^{\theta}}{\lambda_a}\right) \nonumber\\
&&\;\;\leq \;\;\exp\left(n\frac{(\log(2e\lambda_b))^2}{\lambda_a} \right)\exp\left(-\frac{3n}{4} \cdot \frac{(\log{\lambda_b})^{\theta}}{\lambda_a}\right) \nonumber\\
&&\;\;\leq\;\;\exp\left(-\frac{3n}{8} \cdot \frac{(\log{\lambda}_b)^{\theta}}{\lambda_a}\right) \label{chopati}
\end{eqnarray}
for all~\(n \geq N_0\) not depending on~\(\lambda_a\) or~\(\lambda_b,\) where the second inequality in~(\ref{chopati}) is true provided~\(np=\lambda_a \geq \lambda_0\) is large and the final inequality in~(\ref{chopati}) is true if we choose~\(\theta > 2\) strictly.~\(\qed\)

\emph{Proof of Lemma~\ref{thm_inter}~\((b)\)}: Let~\({\cal D}\) be any set containing~\((1+\epsilon)u_n + \Delta\) vertices. Each vertex~\(v \notin {\cal D}\) is adjacent to at least~\((1+\epsilon)u_n\) vertices of~\({\cal D}\) in the graph~\(K_n \setminus H\) and so the probability that~\(v\) is \emph{not} adjacent to any vertex of~\({\cal D}\) in~\(G \setminus H,\) is at most~\((1-p)^{(1+\epsilon)u_n} = \frac{1}{(np)^{1+\epsilon}}.\) Thus if~\({\cal B}\) is the set of all vertices not dominated by~\({\cal D}\) in~\(G \setminus H,\)  then~\(\mathbb{E}\#{\cal B} \leq \frac{n}{(np)^{1+\epsilon}}\) and so a direct application of Markov inequality gives that
\begin{equation}\label{b_est}
\mathbb{P}\left(\#{\cal B} \geq \frac{n \log(np)}{(np)^{1+\epsilon}}\right) \leq \frac{1}{\log(np)}.
\end{equation}
By definition~\({\cal D} \cup {\cal B}\) is a dominating set in~\(G \setminus H\) and has size
\begin{equation}\label{b_est2}
\#({\cal D} \cup {\cal B}) \leq (1+\epsilon)u_n + \Delta + \frac{n \log(np)}{(np)^{1+\epsilon}}
\end{equation}
with probability at least~\(1-\frac{1}{\log{np}},\) by~(\ref{b_est}). Also since~\(p \leq p_0\) a constant, we have that~\(|\log(1-p)|  \leq \sum_{k \geq 1}p^{k} \leq \frac{p}{1-p} \leq \frac{p}{1-p_0}\) and so~\[\frac{n \log(np)}{(np)^{1+\epsilon}} = \frac{1}{(np)^{\epsilon}} \frac{\log(np)}{p} <  \frac{1}{(np)^{\epsilon}} \frac{u_n}{1-p_0} < \epsilon u_n \]
for all~\(n\) large, since~\(np \longrightarrow \infty.\) From~(\ref{b_est2}), we then the upper deviation bound in~(\ref{dom_bound_up_ax}).~\(\qed\)

From the discussion following Lemma~\(1,\) we see that if the edge probability~\(p\) is a constant, then~\(\Delta = o(\log{n})\) is sufficient to ensure the asymptotic equivalence of the robust and the ideal domination numbers. However, in this case we also know that with high probability, the vertex degree in the random graph~\(G\) in fact grows linearly with~\(n,\)  which is much larger than~\(\log{n}.\) Therefore could we, perhaps under additional assumptions, establish asymptotic equivalence for conflict graphs~\(H\) that are sparse \emph{in comparison} to~\(G\)? Addressing this issue, we have the following result for random graphs with a convergent edge probability sequence.
\begin{theorem}\label{thm_main} Suppose~\(np \longrightarrow \infty,\;\;p \leq 1-\frac{1}{n^3}\) and~\(p = p(n) \longrightarrow p_0\) for some constant~\(0\leq p_0 \leq 1.\) As before, let~\(H = H(n)\) be any deterministic graph with maximum vertex degree~\(\Delta = \Delta(n)\) and containing~\(m = m(n)\) edges. If either
\begin{equation}\label{imp_cond}
\Delta=o(n(1-p))\;\;\;\;\text{ or }\;\;\;\;m = o(nu_n(1-p)),
\end{equation}
then~\(\frac{\Gamma_n}{u_n} \longrightarrow 1\) in probability as~\(n \rightarrow \infty.\) 
\end{theorem}
Continuing with constant edge probability example, we see from Theorem~\ref{thm_main} that if~\(p\) is a constant, then either~\(\Delta = o(n)\) or~\(m = o(n\log{n})\) is sufficient for~\(\Gamma_n\) and~\(\gamma(G)\) to be asymptotically equal.

For the case~\(p_0=0\) which is of interest in communication networks, we see that there are robust dominating sets that asymptotically have the same size as the ideal dominating sets even if the number of edges removed per vertex is much larger than the vertex degree itself. This is true because the expected degree of a vertex in~\(G\) equals~\((n-1)p \approx np\) and so by the standard deviation estimate~(\ref{conc_est_f}), we can deduce that each vertex has degree at most of the order of~\(np\) with high probability. Condition~(\ref{imp_cond}) ensures that the robust domination number is asymptotically optimal provided~\(\Delta = o(n),\) even if~\(\Delta\) is much larger than~\(np.\)







To prove Theorem~\ref{thm_main}, we perform a case by case analysis of~\(\Gamma_n\) based on the asymptotic edge probability~\(p_0.\) Recalling the definition of~\(u_n(x,y)\) and~\(u_n = u_n(p,p)\) prior to Lemma~\ref{thm_inter}, we have the following result.
\begin{lemma}\label{thm_dense} Suppose that~\(\Delta \leq r_0 n-1\) for some constant~\(0 < r_0 < 1.\)\\
\((a)\) For every~\(\epsilon > 0\) there are positive constants~\(\lambda_i = \lambda_i(\epsilon,r_0),i=1,2\) such that if~\(\frac{\lambda_1}{n} \leq p \leq \min\left(\frac{1}{2}, 1-\exp\left(-\frac{\epsilon^2(1-r_0)}{16}\right)\right),\) then
\begin{equation}\label{dom_bound_up2}
\mathbb{P}\left(\Gamma_n \leq \frac{(1+6\epsilon)u_n}{1-r_0}\right) \geq 1-z_n,
\end{equation}
where~\[z_n := \min\left(\exp\left(-\frac{\lambda_2 n}{4(np)^{(1+4\epsilon)(1-r_0)^{-1}}}\right), \frac{1}{(np)^{\epsilon/2}}\right).\]
\((b)\) For every~\(\epsilon >0\) and every constant~\(0 < p < 1,\) we have that
\begin{equation}\label{dom_bound_up4}
\mathbb{P}\left(\Gamma_n \leq (1+\epsilon)u_n(p,p(1-\epsilon)-r_0)\right) \geq 1-\exp\left(-\frac{\epsilon^2np}{8}\right).
\end{equation}
\((c)\) For every~\(\epsilon >0,\) there is a constant~\(C = C(\epsilon) > 0\) such that if\\\(q_1:=\max\left(1-p,\frac{\Delta}{n}\right) \leq 2^{-2/\epsilon-2}\) then
\begin{equation}\label{dom_bound_case3}
\mathbb{P}\left(\Gamma_n \leq (1+\epsilon)u_n(p,1-q_1)\right) \geq 1-\frac{C}{n^{\epsilon/2}}.
\end{equation}
\end{lemma}
Parts~\((a),(b)\) and~\((c)\) of Lemma~\ref{thm_dense} essentially obtain deviation upper bounds for~\(\Gamma_n\) for the cases~\(p_0 = 0, 0 < p_0 < 1\) and~\(p_0 =1,\) respectively. 

We use the probabilistic method to prove Lemma~\ref{thm_dense} and so we begin  with a couple of common definitions.  For integer~\(t \geq 1\) let~\({\cal X} := (X_1,\ldots,X_{t})\) be a random~\(t-\)tuple chosen from~\(V^{t},\) that is independent of the graph~\(G.\) Also let~\(\mathbb{P}_X\) denote the probability distribution of~\({\cal X}.\) In each of the three cases below, we choose the tuple~\({\cal X}\) appropriately so that certain niceness properties are satisfied and exploit this to estimate the domination number.

\emph{Proof of Lemma~\ref{thm_dense}~\((a)\)}: For a constant~\(0 < \zeta < 1\) to be determined later, let~\(t = \frac{u_n}{1-\zeta}.\) Assuming that~\(X_i, 1 \leq i \leq t\) are independent and chosen uniformly randomly from~\(V,\) we estimate below the number of vertices ``left out" by the set~\({\cal D} := \{X_1,\ldots,X_t\}.\) For a vertex~\(v,\) the~\(\mathbb{P}_X-\)probability that the random variable~\(X_1\) is equal to~\(v\) or adjacent to~\(v\) in~\(H\) is at most~\(\frac{\Delta+1}{n} \leq r_0.\) Therefore if~\(Q(v)\) is the number of indices~\(i, 1 \leq i \leq t\) such that~\(X_i\) is \emph{not} adjacent to~\(v\) in the graph~\(H,\) then~\(\mathbb{E}_{X}(Q(v)) \geq t(1-r_0)\) and using the standard deviation estimate~(\ref{conc_est_f}) we get for~\(\epsilon > 0\) that
\begin{equation}\label{cruc_eq2}
\mathbb{P}_X\left(Q(v) \geq t(1-r_0)(1-\epsilon)\right) \leq \exp\left(-\frac{\epsilon^2}{4}t(1-r_0)\right).
\end{equation}


In the Appendix we show that~\(u_n(x,x) = \frac{\log(nx)}{|\log(1-x)|}\) is strictly decreasing for all~\(x > \frac{\lambda_{low}}{n}\) where~\(\lambda_{low} > 0\) is a sufficiently large absolute constant. Therefore if~\(\frac{\lambda_{low}}{n} \leq p \leq \lambda_{up} := 1-\exp\left(-\frac{\epsilon^2(1-r_0)}{16}\right),\)  then for all~\(n \geq N_0(\epsilon,r_0)\) large we have that
\begin{equation}\label{t_est_ax}
t \geq u_n(p,p) \geq u_n(\lambda_{up},\lambda_{up}) \geq \frac{\log{n}}{2|\log(1-\lambda_{up})|} = \frac{8}{\epsilon^2(1-r_0)} \log{n},
\end{equation}
where the second inequality in~(\ref{t_est_ax}) is true since~\(\lambda_{up}\) is a constant and so~\(\log(n\lambda_{up}) \geq \frac{\log{n}}{2}\) for all~\(n\) large. Therefore setting~\[E_{tot} := \bigcap_{v \notin {\cal D}} \left\{Q(v) \geq t(1-r_0)(1-\epsilon)\right\},\] we get from the union bound,~(\ref{cruc_eq2}) and~(\ref{t_est_ax}) that~\[\mathbb{P}_X\left(E_{tot}\right) \geq 1-n\cdot \frac{1}{n^2} = 1-\frac{1}{n}.\] We assume henceforth that~\(E_{tot}\) occurs. 


Next, we estimate the number of distinct entries in~\({\cal X}.\) Since~\(|\log(1-p)| > p,\) we get for any constant~\(\zeta > 0\) that~\(t = \frac{1}{1-\zeta}\frac{\log(np)}{|\log(1-p)|} < \frac{\log(np)}{p(1-\zeta)}  < \frac{\epsilon n}{2}(1-r_0),\) provided~\(np \geq \lambda_{low} = \lambda_{low}(\epsilon,\zeta,r_0)\) is large enough. Now, for~\(i \neq j,\) the probability that~\(X_i\) equals~\(X_j\) is~\(\frac{1}{n}\) and so the~\(\mathbb{P}_X-\)expected number of repeated entries is at most~\(\frac{t^2}{n} \leq \frac{\epsilon t}{2}(1-r_0).\) If~\(E_{rep}\) is the event that the number of repeated entries is at most~\(t\epsilon (1-r_0)\) then~\(\mathbb{P}_X(E^c_{rep}) \leq \frac{1}{2},\) by the Markov inequality. Combining with the estimate for~\(\mathbb{P}_X(E_{tot})\)  in the previous paragraph we then get from the union bound that~\(E_{tot} \cap E_{rep}\) occurs with~\(\mathbb{P}_X-\)probability at least~\(\frac{1}{2} - \frac{1}{n} > 0.\)

We assume henceforth that~\(E_{tot} \cap E_{rep}\) occurs so that each vertex~\(v\) is adjacent to least~\(t(1-\epsilon)(1-r_0) - t\epsilon(1-r_0) \geq t(1-2\epsilon) (1-r_0)\) and at most~\(t = \frac{u_n}{1-\zeta}\) vertices  of~\({\cal D},\) in the graph~\( K_n \setminus H.\) Setting~\(1-\zeta := \frac{(1-2\epsilon)(1-r_0)}{1+\epsilon},\) we then get that the probability of the event~\(J_v\) that~\(v\) is not adjacent to any vertex of~\({\cal D}\) in the graph~\(G \setminus H,\) is at most~\((1-p)^{(1+\epsilon)u_n} = \frac{1}{(np)^{1+\epsilon}}\) and at least~\((1-p)^{t} = \frac{1}{(np)^{(1-\zeta)^{-1}}}.\) If~\(L_{tot} := \sum_{v \notin {\cal D}} \ind(J_v),\) then~\(\mathbb{E}L_{tot} \leq \frac{n}{(np)^{1+\epsilon}}\)
and a direct application of the Markov inequality gives us that
\begin{equation}\label{ltot_b}
\mathbb{P}\left(L_{tot} \geq \frac{n}{(np)^{1+\epsilon/2}}\right) \leq \frac{1}{(np)^{\epsilon/2}}.
\end{equation}
Similarly~\(\mathbb{E}L_{tot} \geq \frac{n}{(np)^{(1-\zeta})^{-1}}\) and so using the standard deviation estimate~(\ref{conc_est_f}) we get that
\begin{equation}\label{ltot_a}
\mathbb{P}\left(L_{tot} \geq \frac{2n}{(np)^{1+\epsilon}}\right) \leq \exp\left(-\frac{Cn}{(np)^{(1-\zeta)^{-1}}}\right) \leq \exp\left(-\frac{Cn}{(np)^{(1+4\epsilon)/(1-r_0)}}\right)
\end{equation}
since~\((1-\zeta)^{-1} = \frac{1+\epsilon}{(1-2\epsilon)(1-r_0)} \leq \frac{1+4\epsilon}{1-r_0}\) provided~\(\epsilon > 0\) is a small enough constant. We fix such an~\(\epsilon\) henceforth.

If~\(L_{tot} \leq \frac{n}{(np)^{1+\epsilon/2}},\) then~\(\Gamma_n \leq \#{\cal D} + L_{tot} \leq \frac{u_n}{1-\zeta} + \frac{n}{(np)^{1+\epsilon/2}}\) and moreover, using~\(|\log(1-p)| <2p\) for~\(p < \frac{1}{2},\) we have for~\(np \geq \lambda_{low}\) large enough that~\(\frac{n}{(np)^{1+\epsilon/2}} < \epsilon\frac{\log(np)}{p} < 2\epsilon u_n.\)
Thus~\[\Gamma_n \leq \left(\frac{1}{1-\zeta} + 2\epsilon\right)u_n \leq \left(\frac{1+4\epsilon}{1-r_0} + 2\epsilon\right)u_n \leq \frac{(1+6\epsilon)u_n}{1-r_0}\] and together with~(\ref{ltot_b}) and~(\ref{ltot_a}), this obtains the desired upper bound in~(\ref{dom_bound_up2}).~\(\qed\)


\emph{Proof of Lemma~\ref{thm_dense}~\((b)\)}: We begin with a couple of preliminary calculations. If~\(d_G(v)\) is the degree of~\(v\) in~\(G,\) then~\(\mathbb{E}d_G(v) = (n-1)p.\) Therefore from the deviation estimate~(\ref{conc_est_f}), we get that~\(d_G(v) \geq np(1-\epsilon)\) with probability at least~\(1-\exp\left(-\frac{\epsilon^2}{5}np\right).\) Letting~\(E_{deg} := \bigcap_{v} \{d_G(v) \geq np(1-\epsilon)\},\) we get from the union bound that
\begin{equation}\label{e_deg_est}
\mathbb{P}(E_{deg}) \geq 1-n\exp\left(-\frac{\epsilon^2}{5} np\right) \geq 1- \exp\left(-\frac{\epsilon^2}{8} n p\right)
\end{equation}
for all~\(n\) large.

We henceforth assume that~\(E_{deg}\) occurs and let~\(X_j, 1 \leq j \leq t\) be independently and uniformly chosen from the vertex set~\(V\) also independent of the graph~\(G.\) Let~\({\cal N}(X_i)\) be the set of all neighbours of~\(X_i\) in the graph~\(G\setminus H\) and set~\({\cal N}[X_i] := \{X_i\} \cup {\cal N}(X_i)\) to be the closed neighbourhood of~\(X_i.\) Setting~\({\cal B} := \bigcup_{1 \leq j \leq t} {\cal N}[X_j],\) we see that each vertex in~\({\cal B}\) is adjacent to at least one vertex in~\(\{X_j\}_{1 \leq j \leq t}\) and so~\({\cal D} := \bigcup \{X_j\}_{1 \leq j \leq t} \bigcup \left(V \setminus {\cal B}\right)\)
is a dominating set for~\(G \setminus H.\) Letting~\(\mathbb{P}_X\) be the distribution of~\(\{X_j\}_{1 \leq j \leq t},\) we see that the~\(\mathbb{P}_X-\)expected size of~\({\cal D}\) is
\begin{equation}\label{dom_est}
\mathbb{E}_X\#{\cal D} \leq t + (n-\mathbb{E}_X\#{\cal B})
\end{equation}
and  in the rest of the proof below, we use telescoping to bound the expected size of~\({\cal B}.\)

Formally, for~\(1 \leq j \leq t\) we let~\({\cal B}_j := \bigcup_{1 \leq i \leq j} {\cal N}[X_i]\) and estimate the expected increment~\(\mathbb{E}_X\#{\cal B}_{j} - \mathbb{E}_X\#{\cal B}_{j-1}.\) Adding these increments would then give us the desired bound for~\({\cal B} = {\cal B}_{t}.\) Specifically, we have by construction that~\(\#{\cal B}_j = \#{\cal B}_{j-1} + \#\left({\cal N}[X_j] \setminus {\cal B}_{j-1}\right)\) and for any set~\({\cal S},\)
\begin{eqnarray}
\#\left({\cal N}[X_j] \cap {\cal S}\right) &=& \sum_{y \in {\cal S}} \ind\left(y \in {\cal N}[X_j]\right) \nonumber\\
&=& \sum_{y \in {\cal S}} \ind(y = X_j) + \ind\left(y \in {\cal N}(X_j)\right) \nonumber\\
&=& \sum_{y \in {\cal S}} \ind(y = X_j) + \ind\left(X_j \in {\cal N}(y)\right). \nonumber
\end{eqnarray}
Thus~\(\mathbb{E}_X\#\left({\cal N}[X_j] \cap {\cal S}\right)  = \frac{1}{n} \sum_{y \in {\cal S}} (d(y)+1),\) where~\(d(y)\) is the degree of vertex~\(y\) in~\(G \setminus H\) and setting~\({\cal S} = V \setminus {\cal B}_{j-1} =: {\cal B}^c_{j-1},\) we therefore get that
\begin{equation}\label{somani_one}
\mathbb{E}_X\#{\cal B}_j = \mathbb{E}_X\#{\cal B}_{j-1} + \frac{1}{n}\mathbb{E}\sum_{y \in {\cal B}^c_{j-1}} (d(y)+1).
\end{equation}

We recall that~\(E_{deg}\) occurs and also that the maximum vertex degree of~\(H\) is~\(\Delta\) and so~\(\sum_{y \in {\cal B}^c_{j-1}} (d(y)+1) \geq \#{\cal B}^c_{j-1} (np(1-\epsilon)+1-\Delta).\) Defining~\(\beta_{j} := \frac{\mathbb{E}\#{\cal B}_j}{n},\) we then get that~\(\beta_{j} \geq  \theta_1 + \theta_2 \beta_{j-1},\) where~\(\theta_1 = 1-\theta_2 :=p -p\epsilon -\frac{\Delta-1}{n}.\) Applying recursion and using the fact that~\(\beta_1 > 0,\) we get that
\begin{eqnarray}
\beta_t &\geq& \theta_1\left(1+ \theta_2 + \ldots +\theta_2^{t-2}\right) + \theta_2^{t-1} \beta_1  \nonumber\\
&\geq& \frac{\theta_1}{1-\theta_2}(1-\theta_2^{t-1}) \nonumber\\
&=& 1-(1-\theta_1)^{t-1}.\label{somani_4}
\end{eqnarray}
Substituting~(\ref{somani_4}) into~(\ref{dom_est}) and using~\(\frac{\#{\cal B}}{n} = \beta_{t},\) we finally get that~\(\mathbb{E}_X\#{\cal D} \leq t+ n(1-\theta_1)^{t-1}.\)

Summarizing, if the event~\(E_{deg}\) occurs, then there exists a dominating set of size at most~\(t+ n(1-\theta_1)^{t-1}.\) Setting~\(t -1= (1+\epsilon)\frac{\log{(np)}}{|\log(1-\theta_1)|}\) we get that~\(n(1-\theta_1)^{t-1} =\frac{n}{(np)^{1+\epsilon}} = O\left(\frac{1}{n^{\epsilon}}\right)\) and so the above iteration procedure necessarily terminates after at most~\(t\) steps to provide the desired dominating  set~\({\cal D}\) of size at most~\(t+1.\) By Lemma statement~\[\theta_1 = p(1-\epsilon) - \frac{\Delta-1}{n} \geq p(1-\epsilon)-r_0\] for all~\(n\) large  and so~\({\cal D}\) has size at most~\((1+\epsilon) u_n(p,p(1-\epsilon)-r_0).\) From~(\ref{e_deg_est}), we then get~(\ref{dom_bound_up4}).~\(\qed\)


\emph{Proof of Lemma~\ref{thm_dense}~\((c)\)}: For~\(t \leq 2\log{n},\) let~\({\cal X} = (X_1,\ldots,X_t)\) be a uniformly randomly chosen~\(t-\)tuple from~\(V^{t}\) with \emph{distinct} entries. Say that a vertex~\(v\) is \emph{bad} if~\(v\) is not adjacent to any vertex of~\({\cal X}\) in~\(G \setminus H\) and let~\(q := 1-p.\) The vertex~\(v\) is not adjacent to~\(X_i\) in~\(G \setminus H\) if either~\(v=X_i\) or~\(v\) is adjacent to~\(X_i\) in~\(H\) or the edge~\((v,X_i)\) is not present in~\(G.\) Therefore if~\(A_i\) denotes the event that~\(v\) is not adjacent to~\(X_i\) in~\(G \setminus H,\) then
\begin{equation}\label{pai}
\mathbb{P}(A_i) = \ind({\cal Z}_i) + \ind({\cal Z}_i^c) q
\end{equation}
where~\({\cal Z}_i\) is the event that either~\(v = X_i\) or~\(v\) is adjacent to~\(X_i\) in~\(H.\) Moreover, because the entries of~\({\cal X}\) are distinct, the events~\(A_i\) and~\(A_j\) are mutually~\(\mathbb{P}-\)independent, given~\({\cal X}.\)

Thus denoting~\(A_{bad}(v)\) to be the event that~\(v\) is bad, we see that~
\begin{eqnarray}
\mathbb{E}_X\mathbb{P}(A_{bad}(v)) &=& \mathbb{E}_X\mathbb{P}\left(\bigcap_{1 \leq j \leq t}A_j\right) \nonumber\\
&=& \mathbb{E}_X \prod_{j=1}^{t}\mathbb{P}(A_j) \nonumber\\
 &=& \mathbb{E}_X \left(\prod_{j=1}^{t-1}\mathbb{P}(A_j)\mathbb{E}_X \left(\mathbb{P}(A_t) \mid X_1,\ldots,X_{t-1}\right) \right). \label{neha_tits}
\end{eqnarray}
Given~\(X_1,\ldots,X_{t-1},\) the random variable~\(X_t\) is equally likely to be any of the remaining~\(n-t+1\) vertices from~\(\{1,2,\ldots n\}\) and since the vertex~\(v\) is adjacent to at most~\(\Delta\) vertices in~\(H,\) the event~\({\cal Z}_i\) defined prior to~(\ref{pai}) occurs with conditional probability~\[\mathbb{P}_X({\cal Z}_i \mid X_1,\ldots,X_{t-1}) \leq \frac{\Delta+1}{n-t+1}.\] Plugging this into~(\ref{pai}) we obtain
\[\mathbb{E}_X \left(\mathbb{P}(A_t) \mid X_1,\ldots,X_{t-1}\right) \leq \frac{\Delta+1}{n-t+1} + q \leq \frac{\Delta}{n-2\log{n}} + q\leq 2q_1\] for all~\(n\) large, where~\(q_1 := \max\left(q, \frac{\Delta}{n}\right).\)  Continuing iteratively, we see from~(\ref{neha_tits}) that
\begin{equation}\label{thlop}
\mathbb{E}_X\mathbb{P}(A_{bad}(v)) \leq (2q_1)^{t} \leq \frac{1}{(np)^{1+\epsilon/2}}
\end{equation}
provided we set~\(t := \left(1+\frac{\epsilon}{2}\right) \frac{\log{(np)}}{|\log(2q_1)|}.\) Using~\(p \leq 1\) and~\(q_1 \leq 2^{-2/\epsilon-2}\) we see that the required condition~\(t\leq 2\log{n}\)  is  satisfied for all~\(n\) large.

If~\(N_{bad} := \sum_{v} \ind(A_{bad}(v))\) is the total number of bad vertices, then from~(\ref{thlop}) we see that~\[\mathbb{E}_X \mathbb{E}N_{bad} \leq \frac{n}{(np)^{1+\epsilon/2}}\] and so there exists a choice of~\({\cal X}\) such that~\(\mathbb{E}N_{bad} \leq \frac{n}{(np)^{1+\epsilon/2}}.\) We fix such a~\({\cal X}\) henceforth and get from the Markov inequality that~\[\mathbb{P}(N_{bad} \geq 1) \leq \frac{n}{(np)^{1+\epsilon/2}} \leq \frac{C}{n^{\epsilon/2}},\] for some constant~\(C > 0\) since~\(p \geq 1-2^{-2/\epsilon-2}\) (see statement of the Lemma). In other words, with probability at least~\(1-\frac{C}{n^{\epsilon/2}},\) the vertices in~\({\cal X}\) form a dominating set of~\(G\) and so
\begin{equation}\label{dom_bound_ax_max}
\mathbb{P}\left(\Gamma_n \leq \left(1+\frac{\epsilon}{2}\right) \frac{\log{(np)}}{|\log(2q_1)|}\right) \geq 1-\frac{C}{n^{\epsilon/2}}.
\end{equation}
Again using~\(q_1 \leq 2^{-2/\epsilon-2}\) we have that~\(\frac{1+\epsilon/2}{|\log(2q_1)|} \leq \frac{1+\epsilon}{|\log{q_1}|} \)
and so~(\ref{dom_bound_ax_max}) implies that~\(\Gamma_n \leq (1+\epsilon)u_n(p,1-q_1)\) and this obtains the desired deviation upper bound in~(\ref{dom_bound_up4}).~\(\qed\)

We now use Lemma~\ref{thm_dense} to prove Theorem~\ref{thm_main} below.\\
\emph{Proof of Theorem~\ref{thm_main}}: We consider three separate subcases depending on  whether the asymptotic edge probability~\(p_0 = 0,1\) or otherwise. For~\(p_0 = 0\) and~\(\Delta = o(n),\) we use the lower deviation bound in part~\((a)\) of Lemma~\ref{thm_inter} and the upper deviation bound in part~\((a)\) of Lemma~\ref{thm_dense} to get that~\(\frac{\Gamma_n}{u_n} \longrightarrow 1\) in probability.  Similarly, the cases~\(0 < p_0 < 1\) and~\(p_0=1\) are obtained using parts~\((b)\) and~\((c),\) respectively, of Lemma~\ref{thm_dense}.

For~\(m = o(nu_n(1-p)),\) we include a small ``preprocessing" step. First consider the case~\(p_0 = 0.\) For~\(\epsilon > 0\) let~\({\cal Q}\) be the set of all vertices with degree at most~\(\epsilon n.\) In the proof of Lemma~\ref{thm_dense}\((a),\) we now choose~\(X_i, 1 \leq i \leq t\) uniformly and independently from~\({\cal Q}\) and estimate the number of vertices covered by the set~\({\cal D} := \{X_1,\ldots,X_t\} \cup {\cal Q}^c.\) For~\(\epsilon > 0\) and a vertex~\(v,\) the~\(\mathbb{P}_X-\)probability that~\(X_1\) is equal or adjacent to~\(v\) in~\(H\) is at most~\(\frac{\Delta+1}{\#{\cal Q}} \leq \frac{\epsilon}{1-\epsilon} < 2\epsilon\) since by definition, the set~\({\cal Q}^c\) has size~\(o(u_n) <\epsilon u_n < \epsilon n\) for all~\(n\) large.

If~\(L_{tot}\) is the number of vertices ``left out" by~\(\{X_1,\ldots,X_t\},\) then arguing as in the proof of Lemma~\ref{thm_dense}\((a)\) we get that both~(\ref{ltot_b}) and~(\ref{ltot_a}) holds and so~\[\Gamma_n \leq \#{\cal D} + L_{tot} \leq (1+C_1\epsilon)u_n + \#{\cal Q}^c \leq (1+C_2\epsilon)u_n\] for some constants~\(C_1,C_2>0.\) Since~\(\epsilon > 0\) is arbitrary, we argue  as in the first paragraph of this proof to then get that~\(\mathbb{E}\Gamma_n = u_n(1+o(1)).\) An analogous analysis holds for the cases~\(0 < p_0 <1\) and~\(p_0 = 1\) as well.~\(\qed\)

\renewcommand{\theequation}{\arabic{section}.\arabic{equation}}
\setcounter{equation}{0}
\section{The Sparse Regime}\label{sec_sparse}
In this section, we discuss robust domination in the sparse regime when~\(np \longrightarrow \lambda < \infty.\) We consider the cases~\(\lambda = 0\) and~\(0 < \lambda < \infty\) separately and have the following result regarding the robust domination number.
\begin{theorem}\label{thm_spar}  We have:\\
\((a)\) If~\(np \longrightarrow 0,n^2p \longrightarrow \infty\) and either~\(\Delta = o(n)\) or~\(m = o(n^3p),\) then~\(\frac{4\Gamma_n}{n^2p} \longrightarrow 1\) in probability as~\(n \rightarrow \infty.\)\\
\((b)\) Suppose~\(np \longrightarrow \lambda\) for some~\(0 < \lambda < \infty\) and either~\(\Delta = o(n)\) or~\(m = o(n^2).\) For every~\(\epsilon > 0,\) we have
\begin{equation}\label{inter_gam}
\mathbb{P}\left(a(\lambda)(1-\epsilon) \leq \frac{\gamma(G)}{n} \leq \frac{\Gamma_n}{n} \leq b(\lambda) (1+\epsilon)\right) \longrightarrow 1
\end{equation}
where
\begin{eqnarray}\label{ab_def}
&&a(\lambda) := \left\{
\begin{array}{cc}
\lambda e^{-2\lambda} &\;\;\;\lambda \leq \lambda_0 \\
\;\;\\
\frac{\log{\lambda}-3\log\log{\lambda}}{\lambda}, &\;\;\;\lambda > \lambda_0
\end{array}
\right.
\;\;,\;\;
b(\lambda) := \left\{
\begin{array}{cc}
\frac{\lambda}{4}, &\;\;\;\lambda \leq 1 \\
\;\;\\
\frac{\log{\lambda}+1}{\lambda}, &\;\;\;\lambda > 1
\end{array}
\right.
\;\;,\;\;
\nonumber
\end{eqnarray}
and~\(\lambda_0 > 0\) is an absolute constant not depending on the choice of~\(\lambda\) or~\(H.\)
\end{theorem}
Essentially, for~\(\lambda = 0\) we see that~\(\Gamma_n\) is of the order of~\(n^2p\) while for the ``intermediate" regime~\(0 < \lambda < \infty,\) the robust domination number is of the order of~\(n,\) with high probability.


\emph{Proof of Theorem~\ref{thm_spar}\((a)\)}: If~\(Y_{tot}\) and~\(Z_{tot}\) denote, respectively, the number of edges and the number of isolated edges of~\(G \setminus H,\) then~\(\frac{Y_{tot}}{2} \leq \Gamma_n \leq \frac{Z_{tot}}{2}\) and so it suffices to bound~\(Y_{tot}\) and~\(Z_{tot}.\) The expected number of edges in~\(G\) is~\({n \choose 2}p = \frac{n^2p}{2}(1+o(1))\) and so from the deviation estimate~(\ref{conc_est_f}) in Appendix, we get that
\begin{equation}\label{neha_up}
\mathbb{P}\left(2\Gamma_n \geq Z_{tot} \geq \frac{n^2p}{2}(1+\epsilon)\right) \leq \exp\left(-\frac{\epsilon^2 n^2p}{8}\right).
\end{equation}
This provides an upper bound for~\(\Gamma_n.\)


We now obtain general lower bounds for~\(\Gamma_n\) assuming that~\(np \longrightarrow \lambda\) and~\(\Delta \leq r_0n-1\) for some finite constants~\(0 < r_0 < 1\) and~\(0 \leq \lambda < \infty.\) As discussed before, it suffices to obtain a deviation bound for~\(Y_{tot}\) and we use the second moment method. For an edge~\(e \in K_n \setminus H,\) let~\(A_e\) be the event that~\(e\) is isolated so that~\(Y_{tot} = \sum_{e \in K_n \setminus H} \ind(A_e),\) where~\(\ind(.)\) is the indicator function. Since each vertex has degree at most~\(n,\) we have that
\begin{eqnarray}
\mathbb{P}(A_e) &\geq& p(1-p)^{2n-4} \nonumber\\
&=& \frac{p}{(1-p)^{4}} (1-p)^{2n} \nonumber\\
&=& \frac{p}{(1-p)^4}e^{-2\lambda}(1+o(1)) \nonumber\\
&=& pe^{-2\lambda}(1+o(1)) \label{pae}
\end{eqnarray}
and since~\(\Delta \leq r_0 n\) we have that the number of edges in~\(H\) is~\(m \leq \frac{1}{2}\Delta n \leq \frac{1}{2}r_0n^2.\) Therefore from~(\ref{pae}), we get that
\begin{equation}\label{ey_low}
\mathbb{E}Y_{tot} \geq \left({n \choose 2} - m\right) p e^{-2\lambda}(1+o(1)) \geq \frac{n^2p}{2}e^{-2\lambda}(1-r_0-\epsilon)
\end{equation}
for all~\(n\) large.

Next, the minimum vertex degree in~\(K_n \setminus H\) is~\(n-1-\Delta \geq n(1-r_0)\) and so for distinct edges~\(e_1 \neq e_2\) in~\(K_n \setminus H,\) we argue as in~(\ref{pae}) to get that
\[\mathbb{P}\left(A_{e_1} \cap A_{e_2}\right) \leq p^2(1-p)^{4(n-r_0n-3)}  = p^2 e^{-4\lambda(1-r_0)}(1+o(1)).\] Thus again using~(\ref{pae}), we get that~\(\mathbb{P}\left(A_{e_1} \cap A_{e_2}\right) -\mathbb{P}(A_{e_1})\mathbb{P}(A_{e_2})\) is bounded above by
\begin{equation}
p^2e^{-4\lambda}\left(e^{4\lambda r_0} - 1\right)(1+o(1)) \leq \mathbb{P}(A_{e_1})\mathbb{P}(A_{e_2})\left(e^{4\lambda r_0} - 1\right)(1+o(1)) \nonumber
\end{equation}
and therefore
\begin{eqnarray}
var(Y_{tot}) &=& \sum_{e \in K_n \setminus H} \mathbb{P}(A_e) - \mathbb{P}^2(A_e) + \sum_{e_1 \neq e_2} \mathbb{P}\left(A_{e_1} \cap A_{e_2}\right) -\mathbb{P}(A_{e_1})\mathbb{P}(A_{e_2}) \nonumber\\
&\leq& \sum_{e \in K_n \setminus H} \mathbb{P}(A_e)  + \sum_{e_1 \neq e_2} \mathbb{P}(A_{e_1}) \mathbb{P}(A_{e_2}) \left(e^{4\lambda r_0} - 1\right)(1+o(1)) \nonumber\\
&\leq&  \mathbb{E}Y_{tot} + \left(\mathbb{E}Y_{tot}\right)^2 \left(e^{4\lambda r_0} - 1\right)(1+o(1)). \label{var_up}
\end{eqnarray}

Using~(\ref{var_up}),~(\ref{ey_low}) and the Chebychev inequality, we get for~\(\epsilon > 0\) that
\begin{eqnarray}
\mathbb{P}\left(Y_{tot} \leq \mathbb{E}Y_{tot}(1-\epsilon)\right) &\leq& \frac{var(Y_{tot})}{\epsilon^2 (\mathbb{E}Y_{tot})^2} \nonumber\\
&\leq& \frac{1}{\mathbb{E}Y_{tot}} + \left(e^{4\lambda r_0} - 1\right)(1+o(1)) \nonumber\\
&\leq& \frac{C}{n^2p} + \left(e^{4\lambda r_0} - 1\right)(1+o(1)) \nonumber\\
&\leq& \frac{C}{n^2p} + 2\left(e^{4\lambda r_0} - 1\right)\label{ey_low2}
\end{eqnarray}
where~\(C = C(\lambda,r_0,\epsilon) > 0\) is a constant. Again using~(\ref{ey_low}) and~(\ref{ey_low2}), we get that
\begin{equation}\label{spar_est_ax_low}
\mathbb{P}\left(2\Gamma_n \geq Y_{tot} \geq \frac{n^2p}{2}e^{-2\lambda}(1-r_0-2\epsilon)\right) \geq 1- \frac{C}{n^2p} - 2\left(e^{4\lambda r_0} - 1\right).
\end{equation}

If~\(\Delta = o(n)\) and~\(\lambda = 0,\) then we can set~\(r_0\) arbitrarily small in the above analysis and get from~(\ref{neha_up}) and~(\ref{spar_est_ax_low}) that~\(\frac{4\Gamma_n}{n^2p} \longrightarrow 1\) in probability as~\(n \rightarrow \infty.\) If~\(m = o(n^3p),\) then the number of vertices with degree larger than~\(\epsilon n\) for~\(\epsilon > 0\) is~\(o(n^2p).\) Performing the ``pre-processing" steps as in the proof of Theorem~\ref{thm_main}, we again get that~\(\frac{4\Gamma_n}{n^2p} \longrightarrow 1\) in probability.~\(\qed\)


\emph{Proof of Theorem~\ref{thm_spar}\((b)\)}:  We show that there exists a constant~\(C > 0\) such that for every~\(\epsilon > 0,\)
\begin{equation}\label{egg}
a(\lambda)(1-\epsilon) \leq \frac{\mathbb{E}\Gamma_n}{n} \leq b(\lambda)(1+\epsilon) \text{ and } var(\Gamma_n) \leq C n (\log{n})^2 = o(\mathbb{E}\Gamma_n)^2.
\end{equation}
From~(\ref{egg}) and the Chebychev inequality, we then get~(\ref{inter_gam}). Also, we only consider the case~\(\Delta= o(n)\) and the preprocessing arguments analogous to the proof of Theorem~\ref{thm_spar}\((a)\) holds for the case~\(m = o(n^2).\)

For convenience, we assume throughout that~\(np = \lambda\) and begin with the lower bounds for~\(\mathbb{E}\Gamma_n.\) Set~\(\theta = 3\) in Lemma~\ref{thm_inter}\((a)\) and let~\(\lambda_0 := \lambda_0(3).\) Since~\(p = \frac{\lambda}{n},\) we have that~\[u_n = \frac{\log(np)}{|\log(1-p)|} \sim n\frac{\log{\lambda}}{\lambda}, \lambda_a = np = \lambda, \lambda_b = n|\log(1-p)| \sim \lambda\] and so for~\(\epsilon > 0\) we have that~\(u_n \left(1-\frac{3\log\log{\lambda_b}}{\log{\lambda_a}}\right) \geq a(\lambda)n(1-\epsilon)\) for all~\(n\) large. Consequently, for~\(\lambda > \lambda_0,\) we get from~(\ref{dom_bound_low}) in Lemma~\ref{thm_inter} that~\[\mathbb{P}\left(\Gamma_n \geq a(\lambda)n(1-\epsilon)\right) \geq 1-e^{-Cn}\] and so~\(\mathbb{E}\Gamma_n \geq a(\lambda)n(1-2\epsilon)\) for all~\(n\) large. For~\(\lambda < \lambda_0,\) we use~(\ref{ey_low}) and the fact that~\(\Delta = o(n)\) to get that~\(\mathbb{E}\Gamma_n \geq \frac{n^2p}{4}e^{-2\lambda}(1-2\epsilon) = a(\lambda)n(1-2\epsilon)\) for all~\(n\) large.

Next, for the upper bound for~\(\mathbb{E}\Gamma_n\) for~\(\lambda \leq 1,\) we recall from the discussion prior to~(\ref{ey_low}) that~\(\mathbb{E}\Gamma_n \leq \frac{\mathbb{E}Z_{tot}}{2} \leq \frac{n^2p}{4} = \frac{\lambda n}{4}.\) For~\(\lambda > 1,\) we use the alteration method as in the proof of Lemma~\ref{thm_inter}\((b).\) Let~\({\cal D}\) be any set of~\(u_n + \Delta\) vertices. Each vertex~\(v \notin {\cal D}\) is adjacent to at least~\(u_n\) vertices of~\({\cal D}\) and so~\(v\) is not adjacent to any vertex of~\({\cal D}\) in~\(G \setminus H\) with probability at most~\((1-p)^{u_n} = \frac{1}{np}\) and if~\({\cal B}\) is the set of all ``bad" vertices in~\({\cal D}^c\) not adjacent to any vertex of~\({\cal D},\) then the expected size of~\({\cal B}\) is at most~\(\frac{1}{p} = \frac{n}{\lambda}.\) Moreover, the asymptotic relation~\(|\log(1-p)| \sim p\) and~\(\Delta = o(n)\) imply that~\({\cal D}\) has size~\(u_n + \Delta \leq \left(\frac{\log{\lambda}}{\lambda} + o(1)\right)n.\)  The set~\({\cal D} \cup {\cal B}\) is a dominating set of~\(G \setminus H\) and has an expected size of at most~\(n\left(\frac{\log{\lambda}+1}{\lambda} + o(1)\right) \leq  b(\lambda)n(1+\epsilon)\) for all~\(n\) large. This completes the proof of the expectation bounds in~(\ref{egg}).



Next, to prove the variance bound in~(\ref{egg}), we use the martingale difference method. For~\(1 \leq j \leq n,\) let~\({\cal F}_j = \sigma\left(\{Z(f) : f = (u,v), 1 \leq u < v \leq j\}\right)\) denote the sigma field
generated by the state of the edges in the complete subgraph~\(K_j.\) Defining the martingale difference~\(R_j := \mathbb{E}(\Gamma_n \mid {\cal F}_j) - \mathbb{E}(\Gamma_n\mid {\cal F}_{j-1}),\) we get that~\(\Gamma_n -\mathbb{E}\Gamma_n = \sum_{j=1}^{n} R_j.\) By the martingale property we then have
\begin{equation} \label{var_exp}
var(\Gamma_n)  = \mathbb{E}\left(\sum_{j=1}^{n} R_j\right)^2 = \sum_{j=1}^{n} \mathbb{E}R_j^2.
\end{equation}

To evaluate~\(\mathbb{E}R_j^2,\) we introduce the graph~\(G^{(j)}\) obtained by using independent copies for the states of all edges~\((u,j), 1 \leq u < j\) and retaining the same state as~\(G\) for the rest of the edges. With this notation, we rewrite~\(R_j = \mathbb{E}( \Gamma_n - \Gamma_n^{(j)} \mid {\cal F}_j),\)
where~\(\Gamma_n^{(j)}\)  is the domination number of the graph~\(G^{(j)} \setminus H\) and so squaring and taking expectations, we get~\(\mathbb{E}R_j^2 \leq \mathbb{E}(\Gamma_n - \Gamma_n^{(j)})^2.\) To estimate the difference~\(|\Gamma_n - \Gamma_n^{(j)}|,\) we let~\({\cal D}\) be any minimum size dominating set of~\(G \setminus H.\) Adding all vertices adjacent to the vertex~\(j\) in the graph~\(G^{(j)}\) to the set~\({\cal D}\) gives us a dominating set of~\(G^{(j)} \setminus H\) and so~\(\Gamma_n \leq \Gamma_n^{(j)} + l_j,\) where~\(l_j\) is the total number of edges containing~\(j\) as an endvertex either in~\(G\) or~\(G^{(j)}.\) By symmetry, we therefore get~\(|\Gamma_n - \Gamma^{(j)}_n| \leq l_j\) and so~\(\mathbb{E}R_j^2 \leq \mathbb{E}l^2_j.\)

The expected number of edges in~\(G\) containing~\(j\) as an endvertex is\\\((n-1)p \leq \lambda\) and so if~\(E_{up}\) is the event that the vertex~\(j\) is adjacent to at most~\(C_1 \log{n}\) edges in~\(G\) then using Chernoff bound we get for~\(s > 0\) that
\begin{eqnarray}
\mathbb{P}(E^{c}_{up}) &\leq& e^{-sC_1 \log{n}}(1-p + e^{s}p)^{n-1} \nonumber\\
&\leq& e^{-sC_1 \log{n}}\exp\left((e^{s}-1)(n-1)p\right) \nonumber\\
&\leq& C_2e^{-sC_1 \log{n}} \label{neha_axx}
\end{eqnarray}
for some constant~\(C_2= C_2(\lambda,s) > 0.\) Setting~\(s=1\) and choosing~\(C_1\) large, we get that~\(\mathbb{P}(E^{c}_{up}) \leq \frac{1}{n^{6}}.\) Defining an analogous event~\(E_{up}^{(j)}\) for the graph~\(G^{(j)},\) we get from the union bound that~\(F_{up} := E_{up} \cap E^{(j)}_{up}\) occurs with probability at least~\(1-\frac{2}{n^{6}}.\)  If~\(F_{up}\) occurs, then~\(l_j \leq 2C_1 \log{n}\) and other wise, we use the bound~\(l_j \leq 2n.\) Combining this with the discussion in the previous paragraph, we get
\begin{eqnarray}
\mathbb{E}R_j^2 &\leq& \mathbb{E}l_j^2 \nonumber\\
&\leq& (2C_1\log{n})^2 + 4(2n)^2\mathbb{P}(F_{up}^c)  \nonumber\\
&\leq& (2C_1 \log{n})^2 + \frac{2(2n)^2}{n^{6}} \nonumber\\
&\leq& C_3 (\log{n})^2 \nonumber
\end{eqnarray}
for some constant~\(C_3 > 0.\) Plugging this into~(\ref{var_exp}) gives the variance bound in~(\ref{egg}) and therefore completes the proof of Theorem~\ref{thm_spar}\((b).\)~\(\qed\)

\setcounter{equation}{0}
\renewcommand\theequation{A.\arabic{equation}}
\section*{Appendix}
\emph{Standard Deviation Estimate}: Let~\(Z_i, 1 \leq i \leq t\) be independent Bernoulli random variables satisfying~\(\mathbb{P}(Z_i = 1) = p_i = 1-\mathbb{P}(Z_i = 0).\) If~\(W_t = \sum_{i=1}^{t} Z_i\) and~\(\mu_t = \mathbb{E}W_t,\) then for any~\(0 < \eta < \frac{1}{2}\) we have that
\begin{equation}\label{conc_est_f}
\mathbb{P}\left(\left|W_t-\mu_t\right| \geq \eta \mu_t\right) \leq 2\exp\left(-\frac{\eta^2}{4}\mu_t\right).
\end{equation}
For a proof of~(\ref{conc_est_f}), we refer to Corollary~\(A.1.14,\) pp.~\(312,\) Alon and Spencer (2008).

\underline{\emph{Montonicity of~\(u_n(x)\)}}: The function~\(u_n(x) := \frac{\log(nx)}{|\log(1-x)|}\) has a derivative \[u'_n(x) = \frac{H(x)-x\log{n}}{x(1-x)|\log(1-x)|^2}\] where
\begin{equation}\label{ent_def}
H(x) := - x \cdot \log{x} -(1-x) \cdot \log(1-x)
\end{equation}
is the binary entropy function and logarithms are natural throughout.  If~\(x > \frac{1}{2},\) then~\(H(x) - x\log{n} \leq 1-\frac{\log{n}}{2} <0\) for all~\(n \geq 4.\) The numerator~\(H(x)-x\log{n}\) has  derivative~\(\log\left(\frac{1}{x}-1\right) - \log{n} < 0\) for all~\(x > \frac{1}{n+1}.\) Thus for~\(\frac{\lambda}{n}  < x < \frac{1}{2}\) and~\(\lambda > 1\) we use~\((1-x)|\log(1-x)| <x\) to get that~\(H(x)-x\log{n}\) is bounded above by
\begin{equation}
H\left(\frac{\lambda}{n}\right) - \frac{\lambda\log{n}}{n} = \frac{-\lambda\log{\lambda}}{n} - \left(1-\frac{\lambda}{n}\right)\log\left(1-\frac{\lambda}{n}\right) \leq -\frac{\lambda\log{\lambda}}{n} + \frac{\lambda}{n} \nonumber
\end{equation}
which is strictly less than zero if~\(\lambda > e.\)

\subsection*{\em Data Availability}
Data sharing not applicable to this article as no datasets were generated or analysed during the current study.

\subsection*{\em Acknowledgement}
I thank Professors Rahul Roy, C. R. Subramanian and the referee for crucial comments that led to an improvement of the paper. I also thank IMSc and IISER Bhopal for my fellowships.

\bibliographystyle{plain}

\end{document}